\documentclass[12pt]{article}
\begin{document}
\begin{titlepage}
\begin{center}

{\bf ON THE BOUNDARY VALUE PROBLEM WITH THE OPERATOR IN BOUNDARY CONDITIONS FOR THE OPERATOR-DIFFERENTIAL EQUATION OF THE THIRD ORDER}

\

{\bf Araz R. Aliev, Sevindj F. Babayeva}
\end{center}

\

\centerline{\bf Abstract}

\

{\it In this paper the boundary value problem for one class of the operator-differential equations of the third order on a semi-axis, where one of the boundary conditions is perturbed by some linear operator is researched. There are received sufficient conditions on the operator coefficients of the considered boundary value problem, providing its correct and univalent resolvability in Sobolev type space.}



\end{titlepage}

\newpage


\begin{center}

{\bf ON THE BOUNDARY VALUE PROBLEM WITH THE OPERATOR IN BOUNDARY CONDITIONS FOR THE OPERATOR-DIFFERENTIAL EQUATION OF THE THIRD ORDER}

\

{\bf Araz R. Aliev, Sevindj F. Babayeva}

\end{center}

A number of problems of mathematical physics and mechanics can be reduced to boundary value problems for the differential equations with operators in boundary conditions. In T.Kato's book [1, ch.7] it is possible to meet statements of such problems. In particular, the non-local problem  is one of them. We note, that in works of many mathematicians similar problems for differential equations of the second order are researched in details. Among these works it is possible to specify, for example the works of M.G.Gasymov and S.S.Mirzoev [2], V.A.Ilin and A.F.Filippov [3], M.L.Gorbachuk [4], F.S.Rofe-Beketov [5], S.Y.Yakubov and B.A.Aliev [6], S.S.Mirzoev and Kh.V.Yagubova [7], A.R.Aliev [8]. But we think, that there are few works, devoted in this direction to the equations of the third order which model currents of a liquid in visco-elastic deformable tubes. In this paper we try to fill this gaps. Moreover in comparison with the differential equations of the even order, there are few works, in which the equations of the odd order with scalar boundary conditions on semi-axis are investigated (see, for example, [9-13]).

In the given paper we investigate  the boundary value problem for the operator-differential equation of the third order on semi-axis, where the equation and one of the boundary conditions are perturbed.

{\bf 1.} Let $A$ be a self- adjoined  positive-defined operator in a separable Hilbert space $H$, and $H_{\gamma }$ is a scale of Hilbert spaces, generated by operator $A$, i.e. $D(A^{\gamma})=H_{\gamma},(x,y)_{\gamma}=(A^{\gamma}x,A^{\gamma}y),x,y\in H_{\gamma},(\gamma \ge 0)$. If $\gamma =0,$ we let, that $H_{0} =H$. We denote by $L_{2} \left(\left(a;b\right);H\right),$ $-\infty \le a<b\le +\infty$, a Hilbert space of the vector functions $f(t)$, defined in $\left(a;b\right)$ almost everywhere, with values in $H$, measurable, quadratically integrable in sense of Bochner:
$$
\left\| f\right\| _{L_{2} \left(\left(a;b\right);H\right)} =\, \left(\int _{a}^{b}\left\| f(t)\right\| _{H}^{2}  dt\right)^{{\raise0.5ex\hbox{$\scriptstyle 1 $}\kern-0.1em/\kern-0.15em\lower0.25ex\hbox{$\scriptstyle 2 $}} }.
$$
For $R=(-\infty ;+\infty )$ and $R_{+} =(0;+\infty )$ we assume, that
$$
L_{2} ((-\infty ;+\infty );\, H)\equiv \, L_{2} (R;H),\, \, \, L_{2} ((0;+\infty );\, H)\equiv \, L_{2} (R_{+} ;H).
$$
Further for the vector functions $u(t)$ that almost everywhere belong to $D\left( {A^3 } \right)$ and have the derivative $u'''\left( t \right)$ we determine Hilbert space [14, ch.1]
$$
W_2^3 \left( {R_ + ;H;A} \right)=\left\{u:\, \, u'''\in L_{2} \left(\left(a;b\right);\, H\right)\, ,\, A^{3} u\in L_{2} \left(\left(a;b\right);\, H\right)\right\}
$$
with norm
$$
\left\| u \right\|_{W_2^3 \left( {R_ + ;H;A} \right)}=\left(\left\| u'''\right\| _{L_{2} \left(\left(a;b\right);\, H\right)\, }^{2} +\left\| A^{3} u\right\| _{L_{2} \left(\left(a;b\right);\, H\right)\, }^{2} \right)^{{\raise0.7ex\hbox{$ 1 $}\!\mathord{\left/{\vphantom{1 2}}\right.\kern-\nulldelimiterspace}\!\lower0.7ex\hbox{$ 2 $}} } .
$$
Thus we also accept, that
$$
W_{2}^{3} ((-\infty ;+\infty );H;A)\equiv W_{2}^{3} (R;H;A)\, ,\, \, \, \, W_{2}^{3} ((0,+\infty );H;A)\equiv W_{2}^{3} (R_{+};H;A).
$$
Here all derivatives $u^{(j)} \equiv \frac{d^{j} u}{dt^{j} },\,(j=\overline{1,3})$ are understood in sense of the theory of distributions [14, ch.1].

Let's consider the operators of taking the tracks
$$
\Gamma _{0} u = u(0),\Gamma _{1} u = u'(0), u\in W_{\rm 2}^3 \left( {R_ + ;H;A} \right).
$$
From the theorem of tracks [14, ch.1] it follows, that $\Gamma_{0}:W_{\rm 2}^3 \left( {R_ + ;H;A} \right) \to H_{{\raise0.5ex\hbox{$\scriptstyle 5$}
\kern-0.1em/\kern-0.15em
\lower0.25ex\hbox{$\scriptstyle 2$}}}$, $\Gamma_{1}:W_{\rm 2}^3 \left( {R_ + ;H;A} \right) \to H_{{\raise0.5ex\hbox{$\scriptstyle 3$}
\kern-0.1em/\kern-0.15em
\lower0.25ex\hbox{$\scriptstyle 2$}}}$ are continuous operators. We denote by
$$
\mathop {W_2^3 }\limits^o \left( {R_ + ;H;A} \right) = \left\{ {u:u \in W_2^3 \left( {R_ + ;H;A} \right),\Gamma_{0} u = u(0) = 0,\Gamma_{1} u = u'(0) = 0} \right\}.
$$

Let $L\left(X,Y\right)$ be a space of the bounded operators, acting  from space $X$ to space $Y$.

We also assume, that an operator $K \in L\left( {W_2^3 \left( {R_ + ;H;A} \right),H_{{\raise0.5ex\hbox{$\scriptstyle 3$}
\kern-0.1em/\kern-0.15em \lower0.25ex\hbox{$\scriptstyle 2$}}} } \right)$ and we denote by
$$
\mathop {W_{2;K}^3 }\limits^o \left( {R_ + ;H;A} \right) = \left\{ {u:u \in W_2^3 \left( {R_ + ;H;A} \right),\Gamma_{0} u = u(0) = 0,\Gamma_{1} u = u'(0)= Ku} \right\}.
$$
Obviously, as $\mathop {W_2^3 }\limits^o \left( {R_ + ;H;A} \right)$ and $\mathop {W_{2;K}^3 }\limits^o \left( {R_ + ;H;A} \right)$ are complete subspaces of $W_2^3 \left( {R_ + ;H;A} \right)$.

Now we consider in space $H$ the boundary value problem
\begin{equation} \label{GrindEQ__1_}
u'''(t)-A^{3} u(t)+\sum _{j=1}^{3}A_{j} u^{(3-j)}  (t)=f(t),\, \, \, \, t\in R_{+} ,
\end{equation}
\begin{equation} \label{GrindEQ__2_}
u\left(0\right)=0,\, \, \, \, \, u'\left(0\right)-Ku=0,
\end{equation}
where $f(t)\in L_{2} (R_{+} ;H),\, \, u(t)\in W_{2}^{3} (R_{+} ;H;A),A_{j} \, ,\, j=\overline{1,3},$ are linear, in general, unbounded operators, moreover $A$ is the self-adjoined  positive-defined operator, and the operator $K \in L\left( {W_2^3 \left( {R_ + ;H;A} \right),H_{{\raise0.5ex\hbox{$\scriptstyle 3$}
\kern-0.1em/\kern-0.15em \lower0.25ex\hbox{$\scriptstyle 2$}}} } \right)$, i.e. $\left\| {Ku} \right\|_{H_{{\raise0.5ex\hbox{$\scriptstyle 3$}
\kern-0.1em/\kern-0.15em
\lower0.25ex\hbox{$\scriptstyle 2$}}} }  \le \kappa \left\| u \right\|_{W_2^3 \left( {R_ + ;H;A} \right)}.$

Directly  from the equation (1) and boundary conditions (2) we can see that the, main part of the equation (1)
$$
P_{0} \left(d/dt\right)u\left(t\right)=u'''\left(t\right)-A^{3} u\left(t\right)
$$
is perturbed,
$$
P_{1} \left(d/dt\right)u\left(t\right)=\sum _{j=1}^{3}A_{j} u^{\left(3-j\right)} \left(t\right) ,
$$
and the second boundary condition from (2)
$$u'(0)=0
$$
is perturbed by some operator:
$$
u'(0) - Ku = 0,K \in L\left( {W_2^3 (R_ + ;H;A),H_{{\raise0.5ex\hbox{$\scriptstyle 3$}
\kern-0.1em/\kern-0.15em
\lower0.25ex\hbox{$\scriptstyle 2$}}} } \right).
$$

{\bf Definition 1.} {\it If the vector function} $u\left( t \right) \in W_2^3 \left( {R_ + ;H;A} \right)$ {\it satisfies the equation (1) almost everywhere in} $R_+$, {\it then we say, that u(t) is a regular solution of the equation (1).}

{\bf Definition 2.} {\it If for any} $f\left(t\right)\in L_{2} \left(R_{+} ;H\right)$ {\it there is a regular solution of the equation (1), which satisfies boundary conditions (2) in sense}
$$
\mathop {\lim }\limits_{t \to 0} \left\| {u\left( t \right)} \right\|_{H_{{\raise0.5ex\hbox{$\scriptstyle 5$}
\kern-0.1em/\kern-0.15em
\lower0.25ex\hbox{$\scriptstyle 2$}}} }  = 0,{\rm   }\mathop {\lim }\limits_{t \to 0} \left\| {u'\left( t \right) - Ku} \right\|_{H_{{\raise0.5ex\hbox{$\scriptstyle 3$}
\kern-0.1em/\kern-0.15em
\lower0.25ex\hbox{$\scriptstyle 2$}}} }  = 0,
$$
{\it and the inequality}
$$
\left\| u \right\|_{W_2^3 (R_ + ;H;A)}  \le const\left\| f \right\|_{L_2 (R_ + ;H)}
$$
{\it is fulfilled, then we say, that the boundary value problem (1), (2) is regularly solvable.}

In this paper we study the conditions on coefficients ${\rm }A,A_j ,j = \overline {1,3},$ of the operator-differential equation (1) and on operator $K$, participating in boundary conditions (2), which provide  regular resolvability of the problem (1), (2). The boundary value problem (1), (2) for $K = 0$ is researched in works [9, 11] in various situations.

{\bf 2.} First of all we investigate the main part of the boundary value problem (1), (2) in $H$:
\begin{equation} \label{GrindEQ__3_}
u'''\left(t\right)-A^{3} u\left(t\right)=f\left(t\right),t\in R_{+} ,
\end{equation}
\begin{equation} \label{GrindEQ__4_}
u\left(0\right)=0,\, \, \, \, \, u'\left(0\right)-Ku=0,
\end{equation}
where $f\left( t \right) \in L_2 \left( {R_ + ;H} \right),\,\,u\left( t \right) \in W_2^3 \left( {R_ + ;H;A} \right)$.

Denoting by
$$
{\rm P}_0 u = P_0 \left( {d/dt} \right)u,\,u \in \mathop {W_{2;K}^3 }\limits^o \left( {R_ + ;H;A} \right),
$$
and using a technique [15], we shall prove some auxiliary statements.

{\bf Lemma 1.} {\it Let} $\alpha  > 0,\beta  \in R.$ {\it Then for} $x \in H_{{\raise0.5ex\hbox{$\scriptstyle 5$}
\kern-0.1em/\kern-0.15em
\lower0.25ex\hbox{$\scriptstyle 2$}}}$ {\it the inequality}
$$
\left\| {A^3 e^{ - \alpha At} \sin \beta At\,\,x} \right\|_{L_2 \left( {R_ +  ;H} \right)}^2  \le \left( {\frac{1}{{4\alpha }} - \frac{\alpha }{{4(\alpha ^2  + \beta ^2 )}}} \right)\left\| x \right\|_{H_{{\raise0.5ex\hbox{$\scriptstyle 5$}
\kern-0.1em/\kern-0.15em
\lower0.25ex\hbox{$\scriptstyle 2$}}} }^2.
$$
{\it takes place.}

{\bf P r o o f.} Let $y = A^{{\raise0.5ex\hbox{$\scriptstyle 5$}
\kern-0.1em/\kern-0.15em\lower0.25ex\hbox{$\scriptstyle 2$}}} x \in H.$ Then
$$
\left\| {A^3 e^{ - \alpha At} \sin \beta At\,\,x} \right\|_{L_2 \left( {R_ +  ;H} \right)}^2  = \left\| {A^{{\raise0.5ex\hbox{$\scriptstyle 1$}
\kern-0.1em/\kern-0.15em
\lower0.25ex\hbox{$\scriptstyle 2$}}} e^{ - \alpha At} \sin \beta At\,\,y} \right\|_{L_2 \left( {R_ +  ;H} \right)}^2 =
$$
\begin{equation} \label{GrindEQ__5_}
= \int\limits_0^{ + \infty } {(A^{{\raise0.5ex\hbox{$\scriptstyle 1$}
\kern-0.1em/\kern-0.15em
\lower0.25ex\hbox{$\scriptstyle 2$}}} e^{ - \alpha At} \sin \beta At\,\,y,A^{{\raise0.5ex\hbox{$\scriptstyle 1$}
\kern-0.1em/\kern-0.15em
\lower0.25ex\hbox{$\scriptstyle 2$}}} e^{ - \alpha At} \sin \beta At\,\,y)dt = } \int\limits_0^{ + \infty } {(Ae^{ - 2\alpha At} \sin ^2 \beta At\,\,y,y)dt.}
\end{equation}
Using a spectral decomposition of the  operator $A$ in equality (5), we have:
$$
\int\limits_0^{ + \infty } {\left( {Ae^{ - 2\alpha At} \sin ^2 \beta At\,\,y,y} \right)dt = \int\limits_0^{ + \infty } {\left( {\int\limits_\mu ^{ + \infty } {\sigma e^{ - 2\sigma \alpha t} \sin ^2 \beta \sigma t\left( {dE_\sigma  y,y} \right)} } \right)} } dt =
$$
$$
= \int\limits_\mu ^{ + \infty } {\sigma \left( {\int\limits_0^{ + \infty } {e^{ - 2\sigma \alpha t} \sin ^2 \beta \sigma tdt} } \right)\left( {dE_\sigma  y,y} \right)}
$$
Applying the formula of the integration by parts, we receive
\begin{equation} \label{GrindEQ__6_}
\int\limits_0^{ + \infty } {e^{ - 2\sigma \alpha t} \sin ^2 \beta \sigma tdt}  = \frac{1}{{4\sigma \alpha }} - \frac{1}{2}\int\limits_0^{ + \infty } {e^{ - 2\sigma \alpha t} \cos 2\beta \sigma t} dt.
\end{equation}
Taking into consideration, that $\int\limits_0^{ + \infty } {e^{ - 2\sigma \alpha t} \cos 2\beta \sigma t} dt = \frac{\alpha }{{2\sigma (\alpha ^2  + \beta ^2 )}},$ from (6) we obtain
\begin{equation} \label{GrindEQ__7_}
\int\limits_0^{ + \infty } {e^{ - 2\sigma \alpha t} \sin ^2 \beta \sigma tdt}  = \frac{1}{{4\sigma \alpha }} - \frac{\alpha }{{4\sigma (\alpha ^2  + \beta ^2 )}}.
\end{equation}
Substituting the value of an integral (7) in expression (5), we have:
$$
\left\| {A^3 e^{ - \alpha At} \sin \beta At\, \,x} \right\|_{L_2 \left( {R_ +  ;H} \right)}^2  = \int\limits_0^{ + \infty } {(Ae^{ - 2\alpha At} \sin ^2 \beta At\, \,y,y)dt = }
$$
$$= \int\limits_\mu ^{ + \infty } {\sigma \left( {\frac{1}{{4\sigma \alpha }} - \frac{\alpha }{{4\sigma (\alpha ^2  + \beta ^2 )}}} \right)\left( {dE_\sigma  y,y} \right) = } \left( {\frac{1}{{4\alpha }} - \frac{\alpha }{{4(\alpha ^2  + \beta ^2 )}}} \right)\left\| y \right\|_H^2  =
$$
$$
= \left( {\frac{1}{{4\alpha }} - \frac{\alpha }{{4(\alpha ^2  + \beta ^2 )}}} \right)\left\| {A^{{\raise0.5ex\hbox{$\scriptstyle 5$}
\kern-0.1em/\kern-0.15em
\lower0.25ex\hbox{$\scriptstyle 2$}}} x} \right\|_H^2  = \left( {\frac{1}{{4\alpha }} - \frac{\alpha }{{4(\alpha ^2  + \beta ^2 )}}} \right)\left\| x \right\|_{H_{{\raise0.5ex\hbox{$\scriptstyle 5$}
\kern-0.1em/\kern-0.15em
\lower0.25ex\hbox{$\scriptstyle 2$}}} }^2,
$$
i.e.
$$
\left\| {A^3 e^{ - \alpha At} \sin \beta At\, \,x} \right\|_{L_2 \left( {R_ +  ;H} \right)}^2  \le \left( {\frac{1}{{4\alpha }} - \frac{\alpha }{{4(\alpha ^2  + \beta ^2 )}}} \right)\left\| x \right\|_{H_{{\raise0.5ex\hbox{$\scriptstyle 5$}
\kern-0.1em/\kern-0.15em
\lower0.25ex\hbox{$\scriptstyle 2$}}} }^2.
$$
The lemma is proved.

{\bf Corollary 1.} {\it Taking} $\alpha  = \frac{1}{2},\beta  = \frac{{\sqrt 3 }}{2}$ {\it in the lemma 1, we obtain estimation}
$$
\left\| {A^3 e^{ - \frac{1}{2}At} \sin \frac{{\sqrt 3 }}{2}At\, \,x} \right\|_{L_2 \left( {R_ +  ;H} \right)}^2  \le \frac{3}{8}\left\| x \right\|_{H_{{\raise0.5ex\hbox{$\scriptstyle 5$}
\kern-0.1em/\kern-0.15em
\lower0.25ex\hbox{$\scriptstyle 2$}}} }^2.
$$

{\bf Lemma 2.} {\it Let} $\kappa  = \left\| K \right\|_{W_2^3 \left( {R_ + ;H;A} \right) \to H_{{\raise0.5ex\hbox{$\scriptstyle 3$}
\kern-0.1em/\kern-0.15em
\lower0.25ex\hbox{$\scriptstyle 2$}}} } < 1.$ {\it Then the equation} ${\rm P}_0 u = 0$ {\it has a unique trivial solution from space} $\mathop {W_{2;K}^3 }\limits^o \left( {R_ + ;H;A} \right).$

{\bf P r o o f.} Let $\omega _1  =  - \frac{1}{2} + \frac{{\sqrt 3 }}{2}i$ and $\omega _2  =  - \frac{1}{2} - \frac{{\sqrt 3 }}{2}i.$ The general solution of the equation $P_0 \left( {d/dt} \right)u\left( t \right) = 0$ from space $W_2^3 (R_ + ;H;A)$ has form [9, 15]
$$
u_0 (t) = e^{\omega _1 At} x_1  + e^{\omega _2 At} x_2 ,\,\,\,\,x_1 ,\,\,x_2  \in H_{{\raise0.5ex\hbox{$\scriptstyle 5$}
\kern-0.1em/\kern-0.15em
\lower0.25ex\hbox{$\scriptstyle 2$}}}.
$$
From the condition $u(0) = 0$ we obtain, that $x_1  =  - x_2.$ From the second boundary condition it follows, that $(\omega _1  - \omega _2 )Ax_1  = K\left( {e^{\omega _1 At}  - e^{\omega _2 At} } \right)x_1.$ From here we find that
$$
x_1  = \frac{1}{{i\sqrt 3 }}A^{ - 1} K\left( {e^{\omega _1 At}  - e^{\omega _2 At} } \right)x_1  \equiv \Phi x_1
$$
and also we have, that
$$
\left\| {\Phi x_1 } \right\|_{H_{{\raise0.5ex\hbox{$\scriptstyle 5$}
\kern-0.1em/\kern-0.15em
\lower0.25ex\hbox{$\scriptstyle 2$}}} }  = \left\| {A^{{\raise0.5ex\hbox{$\scriptstyle 5$}
\kern-0.1em/\kern-0.15em
\lower0.25ex\hbox{$\scriptstyle 2$}}} \frac{1}{{i\sqrt 3 }}\left( {A^{ - 1} K\left( {e^{\omega _1 At}  - e^{\omega _2 At} } \right)x_1 } \right)} \right\|_H  \le
$$
\begin{equation} \label{GrindEQ__8_}
\le \frac{1}{{\sqrt 3 }}\left\| K \right\|_{W_2^3 \left( {R_ + ;H;A} \right) \to H_{{\raise0.5ex\hbox{$\scriptstyle 3$}
\kern-0.1em/\kern-0.15em
\lower0.25ex\hbox{$\scriptstyle 2$}}} } \left\| {e^{\omega _1 At} x_1  - e^{\omega _2 At} x_1 } \right\|_{W_2^3 (R_ + ;H;A)}.
\end{equation}
Applying corollary 1, we receive:
$$
\left\| {e^{\omega _1 At} x_1  - e^{\omega _2 At} x_1 } \right\|_{W_2^3 \left( {R_ + ;H;A} \right)}^2  = \left\| {A^3 \left( {e^{\omega _1 At} x_1  - e^{\omega _2 At} x_1 } \right)} \right\|_{L_2 \left( {R_ +  ;H} \right)}^2  +
$$
$$
+ \left\| {\omega _1^3 A^3 e^{\omega _1 At} x_1  - \omega _2^3 A^3 e^{\omega _2 At} x_1 } \right\|_{L_2 \left( {R_ +  ;H} \right)}^2  = 2\left\| {A^3 \left( {e^{\omega _1 At} x_1  - e^{\omega _2 At} x_1 } \right)} \right\|_{L_2 \left( {R_ +  ;H} \right)}^2  =
$$
$$
= 2\left\| {A^3 \left( {e^{\left( { - \frac{1}{2} + i\frac{{\sqrt 3 }}{2}} \right)At} x_1  - e^{\left( { - \frac{1}{2} - i\frac{{\sqrt 3 }}{2}} \right)At} x_1 } \right)} \right\|_{L_2 \left( {R_ +  ;H} \right)}^2  =
$$
$$
= 2\left\| {A^3 e^{ - \frac{1}{2}At} \left( {e^{\frac{{\sqrt 3 }}{2}iAt} x_1  - e^{ - \frac{{\sqrt 3 }}{2}iAt} x_1 } \right)} \right\|_{L_2 \left( {R_ +  ;H} \right)}^2  =
$$
$$
= 8\left\| {A^3 e^{ - \frac{1}{2}At} \sin \frac{{\sqrt 3 }}{2}At\,\,x_1 } \right\|_{L_2 \left( {R_ +  ;H} \right)}^2  \le 8 \cdot \frac{3}{8}\left\| {x_1 } \right\|_{H_{{\raise0.5ex\hbox{$\scriptstyle 5$}
\kern-0.1em/\kern-0.15em
\lower0.25ex\hbox{$\scriptstyle 2$}}} }^2  = 3\left\| {x_1 } \right\|_{H_{{\raise0.5ex\hbox{$\scriptstyle 5$}
\kern-0.1em/\kern-0.15em
\lower0.25ex\hbox{$\scriptstyle 2$}}} }^2.
$$
From here we have
\begin{equation} \label{GrindEQ__9_}
\left\| {e^{\omega _1 At} x_1  - e^{\omega _2 At} x_1 } \right\|_{W_2^3 \left( {R_ + ;H;A} \right)}  \le \sqrt 3 \left\| {x_1 } \right\|_{H_{{\raise0.5ex\hbox{$\scriptstyle 5$}
\kern-0.1em/\kern-0.15em
\lower0.25ex\hbox{$\scriptstyle 2$}}} }.
\end{equation}
Considering the inequality (9) in the equality (8), we get, that
$$
\left\| {\Phi x_1 } \right\|_{H_{{\raise0.5ex\hbox{$\scriptstyle 5$}
\kern-0.1em/\kern-0.15em
\lower0.25ex\hbox{$\scriptstyle 2$}}} }  \le \frac{\kappa }{{\sqrt 3 }} \cdot \sqrt 3 \left\| {x_1 } \right\|_{H_{{\raise0.5ex\hbox{$\scriptstyle 5$}
\kern-0.1em/\kern-0.15em
\lower0.25ex\hbox{$\scriptstyle 2$}}} }  = \kappa \left\| {x_1 } \right\|_{H_{{\raise0.5ex\hbox{$\scriptstyle 5$}
\kern-0.1em/\kern-0.15em
\lower0.25ex\hbox{$\scriptstyle 2$}}} }.
$$
As $\kappa < 1,$ then the operator $E - \Phi$ is invertible in $H_{{\raise0.5ex\hbox{$\scriptstyle 5$}
\kern-0.1em/\kern-0.15em
\lower0.25ex\hbox{$\scriptstyle 2$}}}$ and, we receive, that $x_1  = 0,$ i.e. $u_0 \left( t \right) = 0.$ The lemma is proved.

Now we pass to the basic results of  the  problem (3), (4).

{\bf Theorem 1.} {\it If} $u \in \mathop {W_{2;K}^3 }\limits^o \left( {R_ + ;H;A} \right)$ {\it and} $\kappa  = \left\| K \right\|_{W_2^3 \left( {R_ + ;H;A} \right) \to H_{{\raise0.5ex\hbox{$\scriptstyle 3$}
\kern-0.1em/\kern-0.15em
\lower0.25ex\hbox{$\scriptstyle 2$}}} }  < 1$ {\it then it takes place the  inequality}
\begin{equation} \label{GrindEQ__10_}
\left\| {{\rm P}_0 u} \right\|_{L_2 (R_ + ;H)}^2  \ge (1 - \kappa )\left\| u \right\|_{W_2^3 (R_ + ;H;A)}^2.
\end{equation}

{\bf P r o o f.} Let $u\left( t \right) \in \mathop {W_{2;K}^3 }\limits^o \left( {R_ + ;H;A} \right).$ Then we have:
$$
\left\| {{\rm P}_0 u} \right\|_{L_2 (R_ + ;H)}^2  = \left\| {\frac{{d^3 u}}{{dt^3 }} - A^3 u} \right\|_{L_2 (R_ + ;H)}^2  =
$$
\begin{equation} \label{GrindEQ__11_}
= \left\| {\frac{{d^3 u}}{{dt^3 }}} \right\|_{L_2 (R_ + ;H)}^2  + \left\| {A^3 u} \right\|_{L_2 (R_ + ;H)}^2  - 2{\mathop{\rm Re}\nolimits} \left( {\frac{{d^3 u}}{{dt^3 }},A^3 u} \right)_{L_2 (R_ + ;H)}.
\end{equation}
Applying the formula of the  integration by parts, we receive
$$
\left( {\frac{{d^3 u}}{{dt^3 }},A^3 u} \right)_{L_2 (R_ + ;H)}  =  - \left( {A^{{\raise0.5ex\hbox{$\scriptstyle 1$}
\kern-0.1em/\kern-0.15em
\lower0.25ex\hbox{$\scriptstyle 2$}}} u''(0),A^{{\raise0.5ex\hbox{$\scriptstyle 5$}
\kern-0.1em/\kern-0.15em
\lower0.25ex\hbox{$\scriptstyle 2$}}} u(0)} \right) + \left( {A^{{\raise0.5ex\hbox{$\scriptstyle 3$}
\kern-0.1em/\kern-0.15em
\lower0.25ex\hbox{$\scriptstyle 2$}}} u'(0),A^{{\raise0.5ex\hbox{$\scriptstyle 3$}
\kern-0.1em/\kern-0.15em
\lower0.25ex\hbox{$\scriptstyle 2$}}} u'(0)} \right) -
$$
$$
- \left( {A^{{\raise0.5ex\hbox{$\scriptstyle 5$}
\kern-0.1em/\kern-0.15em
\lower0.25ex\hbox{$\scriptstyle 2$}}} u(0),A^{{\raise0.5ex\hbox{$\scriptstyle 1$}
\kern-0.1em/\kern-0.15em
\lower0.25ex\hbox{$\scriptstyle 2$}}} u''(0)} \right) - \left( {A^3 u,\frac{{d^3 u}}{{dt^3 }}} \right)_{L_2 (R_ + ;H)},
$$
i.e.
\begin{equation} \label{GrindEQ__12_}
2{\mathop{\rm Re}\nolimits} \left( {\frac{{d^3 u}}{{dt^3 }},A^3 u} \right)_{L_2 (R_ + ;H)}  = \left\| {u'(0)} \right\|_{H_{{\raise0.5ex\hbox{$\scriptstyle 3$}
\kern-0.1em/\kern-0.15em
\lower0.25ex\hbox{$\scriptstyle 2$}}} }^2.
\end{equation}
So, for $\kappa  = \left\| K \right\|_{W_2^3 \left( {R_ + ;H;A} \right) \to H_{{\raise0.5ex\hbox{$\scriptstyle 3$}
\kern-0.1em/\kern-0.15em
\lower0.25ex\hbox{$\scriptstyle 2$}}} }  < 1$ in view of (12) from equality (11) we have:
$$
\left\| {{\rm P}_0 u} \right\|_{L_2 (R_ + ;H)}^2  = \left\| u \right\|_{W_2^3 (R_ + ;H;A)}^2  - \left\| {u'(0)} \right\|_{H_{{\raise0.5ex\hbox{$\scriptstyle 3$}
\kern-0.1em/\kern-0.15em
\lower0.25ex\hbox{$\scriptstyle 2$}}} }^2  =
$$
$$
= \left\| u \right\|_{W_2^3 (R_ + ;H;A)}^2  - \left\| {Ku} \right\|_{H_{{\raise0.5ex\hbox{$\scriptstyle 3$}
\kern-0.1em/\kern-0.15em
\lower0.25ex\hbox{$\scriptstyle 2$}}} }^2  \ge \left( {1 - \kappa } \right)\left\| u \right\|_{W_2^3 (R_ + ;H;A)}^2.
$$
Theorem  is proved.

{\bf Theorem 2.} {\it Let} $A$ {\it is the positive-defined self-adjoint operator in} $H\, \, \,\left( {A = A^ *   \ge \mu _0 E} \right),$ $\kappa  = \left\| K \right\|_{W_2^3 \left( {R_ + ;H;A} \right) \to H_{{\raise0.5ex\hbox{$\scriptstyle 3$}
\kern-0.1em/\kern-0.15em
\lower0.25ex\hbox{$\scriptstyle 2$}}} }  < 1.$ {\it Then the operator} ${\rm P}_0 :\mathop {W_{2;K}^3 }\limits^o \left( {R_ + ;H;A} \right) \to L_2 (R_ +  ;H)$ {\it isomorphicly represents} $\mathop {W_{2;K}^3 }\limits^o \left( {R_ + ;H;A} \right)$ {\it on} $L_2 (R_ + ;H).$

{\bf P r o o f.} From lemma 2 it follows, that $Ker{\rm P}_0  = \left\{ 0 \right\}.$ We shall prove, that for any $f\left( t \right) \in L_2 (R_ + ;H)$ there exists $u\left( t \right) \in \mathop {W_{2;K}^3 }\limits^o \left( {R_ + ;H;A} \right),$ such that ${\rm P}_0 u = f,$ i.e. $im{\rm P}_0  = L_2 (R_ + ;H).$

Let's denote by $f_1 (t) = \left\{ \begin{array}{l}
 f(t),t > 0, \\
 0,t < {\rm 0}{\rm ,} \\
 \end{array} \right.$ and $\mathord{\buildrel{\lower3pt\hbox{$\scriptscriptstyle\frown$}}
\over f} _1 (\xi )$ - Fourier transformation of vector function $f_1 (t) \in L_2 (R;H).$ Then the vector function
$$
u_0 (t) = \frac{1}{{\sqrt {2\pi } }}\int\limits_{ - \infty }^{ + \infty } {\left( { - i\xi ^3 E - A^3 } \right)^{ - 1} \mathord{\buildrel{\lower3pt\hbox{$\scriptscriptstyle\frown$}}
\over f} _1 (\xi )e^{i\xi t} d\xi ,\, \, \,t \in R} ,
$$
satisfies the equation $P_0 (d/dt)u(t) = f(t)$ in $R_ +$ almost everywhere. We shall prove, that $u_0 (t) \in W_2^3 (R;H;A).$ From Plansharel theorem it follows, that it's sufficiently to prove, that $A^3 \mathord{\buildrel{\lower3pt\hbox{$\scriptscriptstyle\frown$}}
\over u} _0 (\xi ),$ $\xi ^3 \mathord{\buildrel{\lower3pt\hbox{$\scriptscriptstyle\frown$}}
\over u} _0 (\xi )\in L_2 (R;H),$ where
$$
\mathord{\buildrel{\lower3pt\hbox{$\scriptscriptstyle\frown$}}
\over u} _0 (\xi ) = \frac{1}{{\sqrt {2\pi } }}\int\limits_{ - \infty }^{ + \infty } {u_0 (t)} e^{ - i\xi t} d\xi.
$$
It is obvious, that
$$\left\| {A^3 \mathord{\buildrel{\lower3pt\hbox{$\scriptscriptstyle\frown$}}
\over u} _0 (\xi )} \right\|_{L_2 (R;H)}^2  = \int\limits_{ - \infty }^{ + \infty } {\left\| {A^3 \mathord{\buildrel{\lower3pt\hbox{$\scriptscriptstyle\frown$}}
\over u} _0 (\xi )} \right\|_H^2 d\xi  = \int\limits_{ - \infty }^{ + \infty } {\left\| {A^3 \left( { - i\xi ^3 E - A^3 } \right)^{ - 1} \mathord{\buildrel{\lower3pt\hbox{$\scriptscriptstyle\frown$}}
\over f} _1 (\xi )} \right\|_H^2 } d\xi  \le }
$$
$$
\le \mathop {\sup }\limits_{\xi  \in R} \left\| {A^3 \left( { - i\xi ^3 E - A^3 } \right)^{ - 1} } \right\|^2 \int\limits_{ - \infty }^{ + \infty } {\left\| {\mathord{\buildrel{\lower3pt\hbox{$\scriptscriptstyle\frown$}}
\over f} _1 (\xi )} \right\|_H^2 d\xi  = } \mathop {\sup }\limits_{\xi  \in R} \left\| {A^3 \left( {i\xi ^3 E + A^3 } \right)^{ - 1} } \right\|^2 \left\| {\mathord{\buildrel{\lower3pt\hbox{$\scriptscriptstyle\frown$}}
\over f} _1 (\xi )} \right\|_{L_2 (R;H)}^2  =
$$
$$
= \mathop {\sup }\limits_{\xi  \in R} \left\| {A^3 \left( {i\xi ^3 E + A^3 } \right)^{ - 1} } \right\|^2 \left\| {f_1 } \right\|_{L_2 (R;H)}^2  = \mathop {\sup }\limits_{\xi  \in R} \left\| {A^3 \left( {i\xi ^3 E + A^3 } \right)^{ - 1} } \right\|^2 \left\| f \right\|_{L_2 (R_ +  ;H)}^2 .
$$
Further, from a spectral decomposition of the operator $A$ it follows, that for any $\xi  \in R$
$$
\left\| {A^3 \left( {i\xi ^3 E + A^3 } \right)^{ - 1} } \right\| = \mathop {\sup }\limits_{\mu  \in \sigma (A)} \left| {\mu ^3 \left( {i\xi ^3  + \mu ^3 } \right)^{ - 1} } \right| \le \mathop {\sup }\limits_{\mu  \ge \mu _0 } \left| {\mu ^3 \left( {\xi ^6  + \mu ^6 } \right)^{{\raise0.5ex\hbox{$\scriptstyle { - 1}$}
\kern-0.1em/\kern-0.15em
\lower0.25ex\hbox{$\scriptstyle 2$}}} } \right| \le 1
$$
and $A^3 \mathord{\buildrel{\lower3pt\hbox{$\scriptscriptstyle\frown$}}
\over u} _0 (\xi ) \in L_2 (R;H).$ It may be similarly proved, that $\xi ^3 \mathord{\buildrel{\lower3pt\hbox{$\scriptscriptstyle\frown$}}
\over u} _0 (\xi ) \in L_2 (R;H).$ Hence $u_0 (t) \in W_2^3 (R;H;A).$

Let's denote by $q(t)$ a narrowing of the vector function $u_0 (t)$ on $\left[ {0; + \infty } \right),$ i.e. $q(t) =$ $\left. {u_0 (t)} \right|_{\left[ {0; + \infty } \right)}.$ It is obvious, that $q(t) \in W_2^3 (R_ + ;H;A).$ Therefore from the theorem of tracks [14, ch.1] $q(0) \in H_{{\raise0.5ex\hbox{$\scriptstyle 5$}
\kern-0.1em/\kern-0.15em
\lower0.25ex\hbox{$\scriptstyle 2$}}},$ $q'(0) \in H_{{\raise0.5ex\hbox{$\scriptstyle 3$}
\kern-0.1em/\kern-0.15em
\lower0.25ex\hbox{$\scriptstyle 2$}}},$ $q''(0) \in H_{{\raise0.5ex\hbox{$\scriptstyle 1$}
\kern-0.1em/\kern-0.15em
\lower0.25ex\hbox{$\scriptstyle 2$}}}.$ The  solution of the equation ${\rm P}_0 u = f$ we shall search in the form of
$$
u(t) = q(t) + e^{\omega _1 At} x_1  + e^{\omega _2 At} x_2,
$$
where $\omega _1  =  - \frac{1}{2} + i\frac{{\sqrt 3 }}{2},\, \,\omega _2  =  - \frac{1}{2} - i\frac{{\sqrt 3 }}{2},$ and $x_1 ,x_2  \in H_{{\raise0.5ex\hbox{$\scriptstyle 5$}
\kern-0.1em/\kern-0.15em
\lower0.25ex\hbox{$\scriptstyle 2$}}}$ are unknown vectors which must be determined. From the condition $u\left( t \right) \in W_{2;K}^3 \left( {R_ + ;H;A} \right)$ it follows, that
$$
\left\{ \begin{array}{l}
 q(0) + x_1  + x_2  = 0, \\
 q'(0) + \omega _1 Ax_1  + \omega _2 Ax_2  - K\left( {q(t) + e^{\omega _1 At} x_1  + e^{\omega _2 At} x_2 } \right) = 0. \\
 \end{array} \right.
$$
From here $(E - \Phi )x_1  = \psi,$ where $\psi  = \frac{1}{{i\sqrt 3 }}\left[ {\omega _2 q(0) - A^{ - 1} q'(0) + A^{ - 1} K\left( {q(t) - q(0)e^{\omega _2 At} } \right)} \right] \in H_{{\raise0.5ex\hbox{$\scriptstyle 5$}
\kern-0.1em/\kern-0.15em
\lower0.25ex\hbox{$\scriptstyle 2$}}}.$ From the condition of the theorem we get $\left\| \Phi  \right\|_{H_{{\raise0.5ex\hbox{$\scriptstyle 5$}
\kern-0.1em/\kern-0.15em
\lower0.25ex\hbox{$\scriptstyle 2$}}}  \to H_{{\raise0.5ex\hbox{$\scriptstyle 5$}
\kern-0.1em/\kern-0.15em
\lower0.25ex\hbox{$\scriptstyle 2$}}} }  < 1,$ so $x_1  = (E - \Phi)^{ - 1} \psi  \in H_{{\raise0.5ex\hbox{$\scriptstyle 5$}
\kern-0.1em/\kern-0.15em
\lower0.25ex\hbox{$\scriptstyle 2$}}}.$ Now we can find $x_2  = - q\left( 0 \right) - (E - \Phi)^{ - 1} \psi  \in H_{{\raise0.5ex\hbox{$\scriptstyle 5$}
\kern-0.1em/\kern-0.15em
\lower0.25ex\hbox{$\scriptstyle 2$}}}.$ Consequently, $u \in \mathop {W_{2;K}^3 }\limits^o \left( {R_ + ;H;A} \right)$ and ${\rm P}_0 u = f.$ And on the other hand,
$$
\left\| {{\rm P}_0 u} \right\|_{L_2 \left( {R_ + ;H} \right)}^2  = \left\| {P_0 \left( {d/dt} \right)u} \right\|_{L_2 \left( {R_ + ;H} \right)}^2  = \left\| {\frac{{d^3 u}}{{dt^3 }} - A^3 u} \right\|_{L_2 \left( {R_ + ;H} \right)}^2  \le 2\left\| u \right\|_{W_2^3 \left( {R_ + ;H;A} \right)}^2.
$$
Therefore from Banach theorem  there is an inverse operator ${\rm P}_0^{ - 1}$ and it is bounded. From here it follows, that $\left\| u \right\|_{W_2^3 (R_ + ;H;A)}  \le const\left\| f \right\|_{L_2 (R_ + ;H)}.$ The theorem is proved.

{\bf 3.} As it becomes clear from the theorem 2, the  norms $\left\| u \right\|_{W_2^3 (R_ + ;H;A)}$ and $\left\| {{\rm P}_0 u} \right\|_{L_2 (R_ + ;H)}$ are equivalent in space $\mathop {W_{2;K}^3 }\limits^o \left( {R_ + ;H;A} \right).$ Therefore it is possible to estimate the norms of operators of intermediate derivatives $A^{3 - j} \frac{{d^j }}{{dt^j }}:\mathop {W_{2;K}^3 }\limits^o \left( {R_ + ;H;A} \right) \to L_2 \left( {R_ + ;H} \right),$ $j = \overline {0,2},$ concerning $\left\| {{\rm P}_0 u} \right\|_{L_2 (R_ + ;H)}.$ We shall note, that methods of solutions of the equations without perturbed boundary conditions in problems  with the perturbed  boundary conditions are actually inapplicable. For example, in work [9] for an estimation the norms of operators of the intermediate derivatives having great value at deriving the conditions of resolvability of boundary value problems, the method of factorization which is inapplicable at research of boundary value problems with nonlocal boundary conditions or with the perturbed  boundary conditions is offered. Here for carrying out of such estimations we shall take advantage, as in the work [10], known inequalities from the analysis with combination of the  inequality (10).

The following theorem is true

{\bf Theorem 3.} {\it Let} $\kappa  = \left\| K \right\|_{W_2^3 \left( {R_ + ;H;A} \right) \to H_{{\raise0.5ex\hbox{$\scriptstyle 3$}
\kern-0.1em/\kern-0.15em
\lower0.25ex\hbox{$\scriptstyle 2$}}} }  < 1.$ {\it Then for any} $u \in \mathop {W_{2;K}^3 }\limits^o \left( {R_ + ;H;A} \right)$ {\it following estimations take place:}
\begin{equation} \label{GrindEQ__13_}
\left\| {A^3 u} \right\|_{L_2 (R_ + ;H)}  \le C_0 (\kappa )\left\| {{\rm P}_0 u} \right\|_{L_2 (R_ + ;H)} ,
\end{equation}
\begin{equation} \label{GrindEQ__14_}
\left\| {A^2 u'} \right\|_{L_2 (R_ + ;H)}  \le C_1 (\kappa )\left\| {{\rm P}_0 u} \right\|_{L_2 (R_ + ;H)} ,
\end{equation}
\begin{equation} \label{GrindEQ__15_}
\left\| {Au''} \right\|_{L_2 (R_ + ;H)}  \le C_2 (\kappa )\left\| {{\rm P}_0 u} \right\|_{L_2 (R_ + ;H)} ,
\end{equation}
where
$$
C_0 (\kappa ) = \left( {1 - \kappa } \right)^{{\raise0.5ex\hbox{$\scriptstyle { - 1}$}
\kern-0.1em/\kern-0.15em
\lower0.25ex\hbox{$\scriptstyle 2$}}} ,\,\,C_1 (\kappa ) = \frac{{2^{{\raise0.5ex\hbox{$\scriptstyle 1$}
\kern-0.1em/\kern-0.15em
\lower0.25ex\hbox{$\scriptstyle 3$}}} }}{{3^{{\raise0.5ex\hbox{$\scriptstyle 1$}
\kern-0.1em/\kern-0.15em
\lower0.25ex\hbox{$\scriptstyle 2$}}} }}\left( {1 + \frac{{3\kappa ^{{\raise0.5ex\hbox{$\scriptstyle 2$}
\kern-0.1em/\kern-0.15em
\lower0.25ex\hbox{$\scriptstyle 3$}}} }}{{2^{{\raise0.5ex\hbox{$\scriptstyle 1$}
\kern-0.1em/\kern-0.15em
\lower0.25ex\hbox{$\scriptstyle 3$}}} }}} \right)^{{\raise0.5ex\hbox{$\scriptstyle 1$}
\kern-0.1em/\kern-0.15em
\lower0.25ex\hbox{$\scriptstyle 2$}}} \left( {1 - \kappa } \right)^{{\raise0.5ex\hbox{$\scriptstyle { - 1}$}
\kern-0.1em/\kern-0.15em
\lower0.25ex\hbox{$\scriptstyle 2$}}},
$$
$$
C_2 (\kappa ) = \frac{{2^{{\raise0.5ex\hbox{$\scriptstyle 1$}
\kern-0.1em/\kern-0.15em
\lower0.25ex\hbox{$\scriptstyle 3$}}} }}{{3^{{\raise0.5ex\hbox{$\scriptstyle 1$}
\kern-0.1em/\kern-0.15em
\lower0.25ex\hbox{$\scriptstyle 2$}}} }} \cdot \frac{{1 + 3^{{\raise0.5ex\hbox{$\scriptstyle 1$}
\kern-0.1em/\kern-0.15em
\lower0.25ex\hbox{$\scriptstyle 2$}}} \kappa ^{{\raise0.5ex\hbox{$\scriptstyle 2$}
\kern-0.1em/\kern-0.15em
\lower0.25ex\hbox{$\scriptstyle 3$}}} }}{{\left( {1 - \kappa } \right)^{{\raise0.5ex\hbox{$\scriptstyle 1$}
\kern-0.1em/\kern-0.15em
\lower0.25ex\hbox{$\scriptstyle 2$}}} }}.
$$

{\bf P r o o f.} The validity of the estimation (13) explicitly follows from the inequality (10).

As $u \in \mathop {W_{2;K}^3 }\limits^o \left( {R_ + ;H;A} \right),$ then by the formula of integration by parts we receive:
$$
\left\| {A^2 u'} \right\|_{L_2 \left( {R_ + ;H} \right)}^2  = \int\limits_0^\infty  {\left( {A^2 u',A^2 u'} \right)} dt =  - \int\limits_0^\infty  {\left( {A^3 u,Au''} \right)} dt =
$$
\begin{equation} \label{GrindEQ__16_}
= - \left( {A^3 u,Au''} \right)_{L_2 (R_ + ;H)}  \le \left\| {A^3 u} \right\|_{L_2 (R_ + ;H)} \left\| {Au''} \right\|_{L_2 (R_ + ;H)} .
\end{equation}
Similarly we have:
$$
\left\| {Au''} \right\|_{L_2 (R_ + ;H)}^2  = \int\limits_0^\infty  {\left( {Au'',Au''} \right)} dt =  - \left( {A^{{\raise0.5ex\hbox{$\scriptstyle 3$}
\kern-0.1em/\kern-0.15em
\lower0.25ex\hbox{$\scriptstyle 2$}}} u'(0),A^{{\raise0.5ex\hbox{$\scriptstyle 1$}
\kern-0.1em/\kern-0.15em
\lower0.25ex\hbox{$\scriptstyle 2$}}} u''(0)} \right) -
$$
$$
- \int\limits_0^\infty  {\left( {A^2 u',u'''} \right)} dt \le \left\| {A^{{\raise0.5ex\hbox{$\scriptstyle 3$}
\kern-0.1em/\kern-0.15em
\lower0.25ex\hbox{$\scriptstyle 2$}}} u'(0)} \right\|_H \left\| {A^{{\raise0.5ex\hbox{$\scriptstyle 1$}
\kern-0.1em/\kern-0.15em
\lower0.25ex\hbox{$\scriptstyle 2$}}} u''(0)} \right\|_H  + \left\| {A^2 u'} \right\|_{L_2 (R_ + ;H)} \left\| {u'''} \right\|_{L_2 (R_ + ;H)}  =
$$
$$
= \left\| {Ku} \right\|_{H_{{\raise0.5ex\hbox{$\scriptstyle 3$}
\kern-0.1em/\kern-0.15em
\lower0.25ex\hbox{$\scriptstyle 2$}}} } \left\| {A^{{\raise0.5ex\hbox{$\scriptstyle 1$}
\kern-0.1em/\kern-0.15em
\lower0.25ex\hbox{$\scriptstyle 2$}}} u''\left( 0 \right)} \right\|_H  + \left\| {A^2 u'} \right\|_{L_2 (R_ + ;H)} \left\| {u'''} \right\|_{L_2 (R_ + ;H)}  \le
$$
\begin{equation} \label{GrindEQ__17_}
\le \kappa \left\| u \right\|_{W_2^3 (R_ + ;H;A)} \left\| {A^{{\raise0.5ex\hbox{$\scriptstyle 1$}
\kern-0.1em/\kern-0.15em
\lower0.25ex\hbox{$\scriptstyle 2$}}} u''(0)} \right\|_H  + \left\| {A^2 u'} \right\|_{L_2 (R_ +  ;H)} \left\| {u'''} \right\|_{L_2 (R_ + ;H)} .
\end{equation}
On the other hand,
$$
\left\| {A^{{\raise0.5ex\hbox{$\scriptstyle 1$}
\kern-0.1em/\kern-0.15em
\lower0.25ex\hbox{$\scriptstyle 2$}}} u''(0)} \right\|_H^2  = 2{\mathop{\rm Re}\nolimits} \int\limits_0^\infty  {\left( {Au'',u'''} \right)} dt = 2{\mathop{\rm Re}\nolimits} \left( {Au'',u'''} \right)_{L_2 (R_ + ;H)} ,
$$
i.e.
\begin{equation} \label{GrindEQ__18_}
\left\| {A^{{\raise0.5ex\hbox{$\scriptstyle 1$}
\kern-0.1em/\kern-0.15em
\lower0.25ex\hbox{$\scriptstyle 2$}}} u''(0)} \right\|_H  \le 2^{{\raise0.5ex\hbox{$\scriptstyle 1$}
\kern-0.1em/\kern-0.15em
\lower0.25ex\hbox{$\scriptstyle 2$}}} \left\| {Au''} \right\|_{L_2 (R_ + ;H)}^{{\raise0.5ex\hbox{$\scriptstyle 1$}
\kern-0.1em/\kern-0.15em
\lower0.25ex\hbox{$\scriptstyle 2$}}} \left\| {u'''} \right\|_{L_2 (R_ + ;H)}^{{\raise0.5ex\hbox{$\scriptstyle 1$}
\kern-0.1em/\kern-0.15em
\lower0.25ex\hbox{$\scriptstyle 2$}}}.
\end{equation}
Considering inequalities (16) and (18) in (17), we receive
$$
\left\| {Au''} \right\|_{L_2 \left( {R_ + ;H} \right)}^2  \le \kappa \left\| u \right\|_{W_2^3 (R_ + ;H;A)} 2^{{\raise0.5ex\hbox{$\scriptstyle 1$}
\kern-0.1em/\kern-0.15em
\lower0.25ex\hbox{$\scriptstyle 2$}}} \left\| {Au''} \right\|_{L_2 (R_ + ;H)}^{{\raise0.5ex\hbox{$\scriptstyle 1$}
\kern-0.1em/\kern-0.15em
\lower0.25ex\hbox{$\scriptstyle 2$}}} \left\| {u'''} \right\|_{L_2 (R_ + ;H)}^{{\raise0.5ex\hbox{$\scriptstyle 1$}
\kern-0.1em/\kern-0.15em
\lower0.25ex\hbox{$\scriptstyle 2$}}}  +
$$
$$
+ \left\| {Au''} \right\|_{L_2 (R_ + ;H)}^{{\raise0.5ex\hbox{$\scriptstyle 1$}
\kern-0.1em/\kern-0.15em
\lower0.25ex\hbox{$\scriptstyle 2$}}} \left\| {A^3 u} \right\|_{L_2 \left( {R_ + ;H} \right)}^{{\raise0.5ex\hbox{$\scriptstyle 1$}
\kern-0.1em/\kern-0.15em
\lower0.25ex\hbox{$\scriptstyle 2$}}} \left\| {u'''} \right\|_{L_2 (R_ + ;H)} ,
$$
i.e.
$$
\left\| {Au''} \right\|_{L_2 (R_ + ;H)}^{{\raise0.5ex\hbox{$\scriptstyle 3$}
\kern-0.1em/\kern-0.15em
\lower0.25ex\hbox{$\scriptstyle 2$}}}  \le 2^{{\raise0.5ex\hbox{$\scriptstyle 1$}
\kern-0.1em/\kern-0.15em
\lower0.25ex\hbox{$\scriptstyle 2$}}} \kappa \left\| u \right\|_{W_2^3 (R_ + ;H;A)} \left\| {u'''} \right\|_{L_2 (R_ + ;H)}^{{\raise0.5ex\hbox{$\scriptstyle 1$}
\kern-0.1em/\kern-0.15em
\lower0.25ex\hbox{$\scriptstyle 2$}}}  + \left\| {A^3 u} \right\|_{L_2 (R_ + ;H)}^{{\raise0.5ex\hbox{$\scriptstyle 1$}
\kern-0.1em/\kern-0.15em
\lower0.25ex\hbox{$\scriptstyle 2$}}} \left\| {u'''} \right\|_{L_2 \left( {R_ + ;H} \right)}.
$$
Taking into consideration, that $\left\| {u'''} \right\|_{L_2 \left( {R_ + ;H} \right)}  \le \left\| u \right\|_{W_2^3 (R_ + ;H;A)},$ we obtain:
\begin{equation} \label{GrindEQ__19_}
\left\| {Au''} \right\|_{L_2 (R_ + ;H)}  \le 2^{{\raise0.5ex\hbox{$\scriptstyle 1$}
\kern-0.1em/\kern-0.15em
\lower0.25ex\hbox{$\scriptstyle 3$}}} \kappa ^{{\raise0.5ex\hbox{$\scriptstyle 2$}
\kern-0.1em/\kern-0.15em
\lower0.25ex\hbox{$\scriptstyle 3$}}} \left\| u \right\|_{W_2^3 (R_ + ;H;A)}  + \left\| {A^3 u} \right\|_{L_2 (R_ + ;H)}^{{\raise0.5ex\hbox{$\scriptstyle 1$}
\kern-0.1em/\kern-0.15em
\lower0.25ex\hbox{$\scriptstyle 3$}}} \left\| {u'''} \right\|_{L_2 (R_ + ;H)}^{{\raise0.5ex\hbox{$\scriptstyle 2$}
\kern-0.1em/\kern-0.15em
\lower0.25ex\hbox{$\scriptstyle 3$}}}.
\end{equation}
And from here
$$
\left( {\left\| {Au''} \right\|_{L_2 (R_ + ;H)}  - 2^{{\raise0.5ex\hbox{$\scriptstyle 1$}
\kern-0.1em/\kern-0.15em
\lower0.25ex\hbox{$\scriptstyle 3$}}} \kappa ^{{\raise0.5ex\hbox{$\scriptstyle 2$}
\kern-0.1em/\kern-0.15em
\lower0.25ex\hbox{$\scriptstyle 3$}}} \left\| u \right\|_{W_2^3 (R_ + ;H;A)} } \right)^2  \le \left\| {A^3 u} \right\|_{L_2 (R_ + ;H)}^{{\raise0.5ex\hbox{$\scriptstyle 2$}
\kern-0.1em/\kern-0.15em
\lower0.25ex\hbox{$\scriptstyle 3$}}} \left\| {u'''} \right\|_{L_2 (R_ + ;H)}^{{\raise0.5ex\hbox{$\scriptstyle 4$}
\kern-0.1em/\kern-0.15em
\lower0.25ex\hbox{$\scriptstyle 3$}}} .
$$
Then for any $\varepsilon > 0$, applying Young inequality we receive:
$$
\left( {\left\| {Au''} \right\|_{L_2 (R_ + ;H)}  - 2^{{\raise0.5ex\hbox{$\scriptstyle 1$}
\kern-0.1em/\kern-0.15em
\lower0.25ex\hbox{$\scriptstyle 3$}}} \kappa ^{{\raise0.5ex\hbox{$\scriptstyle 2$}
\kern-0.1em/\kern-0.15em
\lower0.25ex\hbox{$\scriptstyle 3$}}} \left\| u \right\|_{W_2^3 (R_ + ;H;A)} } \right)^2  \le \left( {\varepsilon \left\| {A^3 u} \right\|_{L_2 \left( {R_ + ;H} \right)}^2 } \right)^{{\raise0.5ex\hbox{$\scriptstyle 1$}
\kern-0.1em/\kern-0.15em
\lower0.25ex\hbox{$\scriptstyle 3$}}} \left( {\frac{1}{{\varepsilon ^{{\raise0.5ex\hbox{$\scriptstyle 1$}
\kern-0.1em/\kern-0.15em
\lower0.25ex\hbox{$\scriptstyle 2$}}} }}\left\| {u'''} \right\|_{L_2 \left( {R_ + ;H} \right)}^2 } \right)^{{\raise0.5ex\hbox{$\scriptstyle 2$}
\kern-0.1em/\kern-0.15em
\lower0.25ex\hbox{$\scriptstyle 3$}}} \le
$$
$$
\le \frac{1}{3}\varepsilon \left\| {A^3 u} \right\|_{L_2 \left( {R_ + ;H} \right)}^2  + \frac{2}{{3\varepsilon ^{{\raise0.5ex\hbox{$\scriptstyle 1$}
\kern-0.1em/\kern-0.15em
\lower0.25ex\hbox{$\scriptstyle 2$}}} }}\left\| {u'''} \right\|_{L_2 (R_ + ;H)}^2 .
$$
Supposing $\frac{1}{3}\varepsilon  = \frac{2}{3}\varepsilon ^{{\raise0.5ex\hbox{$\scriptstyle { - 1}$}
\kern-0.1em/\kern-0.15em
\lower0.25ex\hbox{$\scriptstyle 2$}}} ,$ it we obtain, that $\varepsilon  = 2^{{\raise0.5ex\hbox{$\scriptstyle 2$}
\kern-0.1em/\kern-0.15em
\lower0.25ex\hbox{$\scriptstyle 3$}}}.$ Thus,
$$
\left( {\left\| {Au''} \right\|_{L_2 (R_ + ;H)}  - 2^{{\raise0.5ex\hbox{$\scriptstyle 1$}
\kern-0.1em/\kern-0.15em
\lower0.25ex\hbox{$\scriptstyle 3$}}} \kappa ^{{\raise0.5ex\hbox{$\scriptstyle 2$}
\kern-0.1em/\kern-0.15em
\lower0.25ex\hbox{$\scriptstyle 3$}}} \left\| u \right\|_{W_2^3 (R_ + ;H;A)} } \right)^2  \le \frac{{2^{{\raise0.5ex\hbox{$\scriptstyle 2$}
\kern-0.1em/\kern-0.15em
\lower0.25ex\hbox{$\scriptstyle 3$}}} }}{3}\left\| u \right\|_{W_2^3 (R_ + ;H;A)}^2.
$$
Hence,
$$
\left\| {Au''} \right\|_{L_2 (R_ + ;H)}  \le \left( {\frac{{2^{{\raise0.5ex\hbox{$\scriptstyle 1$}
\kern-0.1em/\kern-0.15em
\lower0.25ex\hbox{$\scriptstyle 3$}}} }}{{3^{{\raise0.5ex\hbox{$\scriptstyle 1$}
\kern-0.1em/\kern-0.15em
\lower0.25ex\hbox{$\scriptstyle 2$}}} }} + 2^{{\raise0.5ex\hbox{$\scriptstyle 1$}
\kern-0.1em/\kern-0.15em
\lower0.25ex\hbox{$\scriptstyle 3$}}} \kappa ^{{\raise0.5ex\hbox{$\scriptstyle 2$}
\kern-0.1em/\kern-0.15em
\lower0.25ex\hbox{$\scriptstyle 3$}}} } \right)\left\| u \right\|_{W_2^3 (R_ + ;H;A)}  =
$$
$$
= \frac{{2^{{\raise0.5ex\hbox{$\scriptstyle 1$}
\kern-0.1em/\kern-0.15em
\lower0.25ex\hbox{$\scriptstyle 3$}}} }}{{3^{{\raise0.5ex\hbox{$\scriptstyle 1$}
\kern-0.1em/\kern-0.15em
\lower0.25ex\hbox{$\scriptstyle 2$}}} }}\left( {1 + 3^{{\raise0.5ex\hbox{$\scriptstyle 1$}
\kern-0.1em/\kern-0.15em
\lower0.25ex\hbox{$\scriptstyle 2$}}} \kappa ^{{\raise0.5ex\hbox{$\scriptstyle 2$}
\kern-0.1em/\kern-0.15em
\lower0.25ex\hbox{$\scriptstyle 3$}}} } \right)\left\| u \right\|_{W_2^3 (R_ + ;H;A)}  \le \frac{{2^{{\raise0.5ex\hbox{$\scriptstyle 1$}
\kern-0.1em/\kern-0.15em
\lower0.25ex\hbox{$\scriptstyle 3$}}} }}{{3^{{\raise0.5ex\hbox{$\scriptstyle 1$}
\kern-0.1em/\kern-0.15em
\lower0.25ex\hbox{$\scriptstyle 2$}}} }} \cdot \frac{{1 + 3^{{\raise0.5ex\hbox{$\scriptstyle 1$}
\kern-0.1em/\kern-0.15em
\lower0.25ex\hbox{$\scriptstyle 2$}}} \kappa ^{{\raise0.5ex\hbox{$\scriptstyle 2$}
\kern-0.1em/\kern-0.15em
\lower0.25ex\hbox{$\scriptstyle 3$}}} }}{{\left( {1 - \kappa } \right)^{{\raise0.5ex\hbox{$\scriptstyle 1$}
\kern-0.1em/\kern-0.15em
\lower0.25ex\hbox{$\scriptstyle 2$}}} }}\left\| {{\rm P}_0 u} \right\|_{L_2 (R_ + ;H)} .
$$
Thus, the estimation (15) is true. Now we shall prove (14). Considering an inequality (19) in (16) and spending the same reasoning, as above, we receive:
$$
\left\| {A^2 u'} \right\|_{L_2 (R_ + ;H)}^2  \le 2^{{\raise0.5ex\hbox{$\scriptstyle 1$}
\kern-0.1em/\kern-0.15em
\lower0.25ex\hbox{$\scriptstyle 3$}}} \kappa ^{{\raise0.5ex\hbox{$\scriptstyle 2$}
\kern-0.1em/\kern-0.15em
\lower0.25ex\hbox{$\scriptstyle 3$}}} \left\| u \right\|_{W_2^3 (R_ + ;H;A)} \left\| {A^3 u} \right\|_{L_2 (R_ + ;H)}  + \left\| {A^3 u} \right\|_{L_2 (R_ + ;H)}^{{\raise0.5ex\hbox{$\scriptstyle 4$}
\kern-0.1em/\kern-0.15em
\lower0.25ex\hbox{$\scriptstyle 3$}}} \left\| {u'''} \right\|_{L_2 (R_ + ;H)}^{{\raise0.5ex\hbox{$\scriptstyle 2$}
\kern-0.1em/\kern-0.15em
\lower0.25ex\hbox{$\scriptstyle 3$}}}  \le
$$
$$
\le 2^{{\raise0.5ex\hbox{$\scriptstyle 1$}
\kern-0.1em/\kern-0.15em
\lower0.25ex\hbox{$\scriptstyle 3$}}} \kappa ^{{\raise0.5ex\hbox{$\scriptstyle 2$}
\kern-0.1em/\kern-0.15em
\lower0.25ex\hbox{$\scriptstyle 3$}}} \left\| u \right\|_{W_2^3 (R_ + ;H;A)}^2  + \left( {\varepsilon \left\| {u'''} \right\|_{L_2 \left( {R_ + ;H} \right)}^2 } \right)^{{\raise0.5ex\hbox{$\scriptstyle 1$}
\kern-0.1em/\kern-0.15em
\lower0.25ex\hbox{$\scriptstyle 3$}}} \left( {\frac{1}{{\varepsilon ^{{\raise0.5ex\hbox{$\scriptstyle 1$}
\kern-0.1em/\kern-0.15em
\lower0.25ex\hbox{$\scriptstyle 2$}}} }}\left\| {A^3 u} \right\|_{L_2 \left( {R_ + ;H} \right)}^2 } \right)^{{\raise0.5ex\hbox{$\scriptstyle 2$}
\kern-0.1em/\kern-0.15em
\lower0.25ex\hbox{$\scriptstyle 3$}}}  \le
$$
$$
\le 2^{{\raise0.5ex\hbox{$\scriptstyle 1$}
\kern-0.1em/\kern-0.15em
\lower0.25ex\hbox{$\scriptstyle 3$}}} \kappa ^{{\raise0.5ex\hbox{$\scriptstyle 2$}
\kern-0.1em/\kern-0.15em
\lower0.25ex\hbox{$\scriptstyle 3$}}} \left\| u \right\|_{W_2^3 (R_ + ;H;A)}^2  + \frac{1}{3}\varepsilon \left\| {u'''} \right\|_{L_2 \left( {R_ + ;H} \right)}^2  + \frac{2}{{3\varepsilon ^{{\raise0.5ex\hbox{$\scriptstyle 1$}
\kern-0.1em/\kern-0.15em
\lower0.25ex\hbox{$\scriptstyle 2$}}} }}\left\| {A^3 u} \right\|_{L_2 \left( {R_ + ;H} \right)}^2 .
$$
Supposing here also $\varepsilon  = 2^{{\raise0.5ex\hbox{$\scriptstyle 2$}
\kern-0.1em/\kern-0.15em
\lower0.25ex\hbox{$\scriptstyle 3$}}},$ we have:
$$
\left\| {A^2 u'} \right\|_{L_2 \left( {R_ + ;H} \right)}^2  \le 2^{{\raise0.5ex\hbox{$\scriptstyle 1$}
\kern-0.1em/\kern-0.15em
\lower0.25ex\hbox{$\scriptstyle 3$}}} \kappa ^{{\raise0.5ex\hbox{$\scriptstyle 2$}
\kern-0.1em/\kern-0.15em
\lower0.25ex\hbox{$\scriptstyle 3$}}} \left\| u \right\|_{W_2^3 (R_ + ;H;A)}^2  + \frac{{2^{{\raise0.5ex\hbox{$\scriptstyle 2$}
\kern-0.1em/\kern-0.15em
\lower0.25ex\hbox{$\scriptstyle 3$}}} }}{3}\left\| u \right\|_{W_2^3 (R_ + ;H;A)}^2  =
$$
$$
= \frac{{2^{{\raise0.5ex\hbox{$\scriptstyle 2$}
\kern-0.1em/\kern-0.15em
\lower0.25ex\hbox{$\scriptstyle 3$}}} }}{3}\left( {1 + \frac{{3\kappa ^{{\raise0.5ex\hbox{$\scriptstyle 2$}
\kern-0.1em/\kern-0.15em
\lower0.25ex\hbox{$\scriptstyle 3$}}} }}{{2^{{\raise0.5ex\hbox{$\scriptstyle 1$}
\kern-0.1em/\kern-0.15em
\lower0.25ex\hbox{$\scriptstyle 3$}}} }}} \right)\left\| u \right\|_{W_2^3 (R_ + ;H;A)}^2  \le \frac{{2^{{\raise0.5ex\hbox{$\scriptstyle 2$}
\kern-0.1em/\kern-0.15em
\lower0.25ex\hbox{$\scriptstyle 3$}}} }}{3}\left( {1 + \frac{{3\kappa ^{{\raise0.5ex\hbox{$\scriptstyle 2$}
\kern-0.1em/\kern-0.15em
\lower0.25ex\hbox{$\scriptstyle 3$}}} }}{{2^{{\raise0.5ex\hbox{$\scriptstyle 1$}
\kern-0.1em/\kern-0.15em
\lower0.25ex\hbox{$\scriptstyle 3$}}} }}} \right)\left( {1 - \kappa } \right)^{ - 1} \left\| {{\rm P}_0 u} \right\|_{L_2 \left( {R_ + ;H} \right)}^2.
$$
So, the estimation (14) is also proved. The theorem is proved.

The  estimations of norms of operators of intermediate derivatives in theorem  3 have also independent mathematical interest. Similar problems for numerical functions can be found and studied, for example, in work [16] and in available there references.

{\bf 4.} Before passing to  establishment of conditions of regular resolvability for the boundary value problem (1), (2), we shall prove the following statement.

{\bf Lemma 3.} {\it Let} $B_j  = A_j A^{ - j},$ $j = \overline {1,3},$ {\it are bounded operators in} $H.$ {\it Then an operator} ${\rm P} = {\rm P}_0  + {\rm P}_1,$ {\it where} ${\rm P}_1$ {\it - is the operator acting by the following way:}
$$
{\rm P}_1 u = P_1 \left( {d/dt} \right)u,\,\,\,u \in \mathop {W_{2;K}^3 }\limits^o \left( {R_ + ;H;A} \right),
$$
{\it is the bounded operator from} $\mathop {W_{2;K}^3 }\limits^o \left( {R_ + ;H;A} \right)$ {\it to} $L_2 (R_ + ;H).$

{\bf P r o o f.} Really, for any $u(t) \in \mathop {W_{2;K}^3 }\limits^o \left( {R_ + ;H;A} \right)$
$$
\left\| {{\rm P}u} \right\|_{L_2 (R_ + ;H)}  \le \left\| {{\rm P}_0 u} \right\|_{L_2 (R_ + ;H)}  + \left\| {{\rm P}_1 u} \right\|_{L_2 (R_ + ;H)}  \le \left\| {{\rm P}_0 u} \right\|_{L_2 (R_ + ;H)}  +
$$
$$
+ \sum\limits_{j = 1}^3 {\left\| {A_j u^{(3 - j)} } \right\|_{L_2 (R_ + ;H)} }  \le \left\| {{\rm P}_0 u} \right\|_{L_2 (R_ + ;H)}  + \sum\limits_{j = 1}^3 {\left\| {B_j } \right\|_{H \to H} \left\| {A^j u^{(3 - j)} } \right\|_{L_2 (R_ + ;H)} } .
$$
Then from this inequality, taking into consideration theorem 2 and theorem of intermediate derivatives [14, ch.1], we receive
$$
\left\| {{\rm P}u} \right\|_{L_2 (R_ + ;H)}  \le const\left\| u \right\|_{W_2^3 (R_ + ;H;A)} .
$$
The lemma is proved.

And now we formulate the basic theorem of regular solvability of the problem (1), (2).

{\bf Theorem 4.} {\it Let the conditions of theorem 2 are satisfied, and operators} $B_j  = A_j A^{ - j},$ $j = \overline {1,3},$ {\it are bounded  in} $H$ {\it and the inequality}
$$
\alpha \left( \kappa  \right) = \sum\limits_{j = 0}^2 {C_j \left( \kappa  \right)\left\| {B_{3 - j} } \right\|_{H \to H}  < 1}
$$
{\it takes place, where} $C_j \left( \kappa  \right),\,\,j = \overline {0,2},$ {\it are defined in theorem 3. Then the problem (1), (2) is regularly solvable.}

{\bf P r o o f.} From the theorem 2 the operator ${\rm P}_0 :\mathop {W_{2;K}^3 }\limits^o \left( {R_ + ;H;A} \right) \to L_2 \left( {R_ + ;H} \right)$ is isomorphism. Then there is a bounded inverse operator ${\rm P}_0^{ - 1}.$ We rewrite the problem (1), (2) in the form of the operator equation ${\rm P}u = {\rm P}_0 u + {\rm P}_1 u = f,$ where $f \in L_2 (R_ + ;H),$ $u \in \mathop {W_{2;K}^3 }\limits^o \left( {R_ + ;H;A} \right).$ After replacement ${\rm P}_0 u = v$ we receive the equation $v + {\rm P}_1 {\rm P}_0^{ - 1} v = f$ from $L_2 (R_ +  ;H).$ But for any $v \in L_2 (R_ +  ;H),$ considering the theorem 3,
$$
\left\| {{\rm P}_1 {\rm P}_0^{ - 1} v} \right\|_{L_2 (R_ + ;H)}  = \left\| {{\rm P}_1 u} \right\|_{L_2 (R_ + ;H)}  \le \left\| {\sum\limits_{j = 0}^2 {A_{3 - j} u^{(j)} } } \right\|_{L_2 (R_ + ;H)}  \le
$$
$$
\le \sum\limits_{j = 0}^2 {\left\| {B_{3 - j} } \right\|} \left\| {A^{3 - j} u^{\left( j \right)} } \right\|_{L_2 (R_ + ;H)}  \le \sum\limits_{j = 0}^2 {C_j \left( \kappa  \right)\left\| {B_{3 - j} } \right\|_{H \to H} }  = \alpha (\kappa ) < 1.
$$
Thus, the operator $E + {\rm P}_1 {\rm P}_0^{ - 1}$ is invertible in $L_2 (R_ +  ;H).$ Then $v = \left( {E + {\rm P}_1 {\rm P}_0^{ - 1} } \right)^{ - 1} f$ and $u = {\rm P}_0^{ - 1} \left( {E + {\rm P}_1 {\rm P}_0^{ - 1} } \right)^{ - 1} f.$ From here it follows, that
$$
\left\| u \right\|_{W_2^3 (R_ + ;H;A)}  \le const\left\| f \right\|_{L_2 (R_ + ;H)}.
$$
The theorem is proved.

{\bf Corollary 2.} {\it Let} $K = 0.$ {\it If the inequality}
$$
\alpha \left( 0 \right) = \frac{{2^{{\raise0.5ex\hbox{$\scriptstyle 1$}
\kern-0.1em/\kern-0.15em
\lower0.25ex\hbox{$\scriptstyle 3$}}} }}{{3^{{\raise0.5ex\hbox{$\scriptstyle 1$}
\kern-0.1em/\kern-0.15em
\lower0.25ex\hbox{$\scriptstyle 2$}}} }}\left( {\left\| {B_1 } \right\|_{H \to H}  + \left\| {B_2 } \right\|_{H \to H} } \right) + \left\| {B_3 } \right\|_{H \to H}  < 1
$$
{\it takes place, the problem (1), (2) is regular solvable.}

We must note, that for $K = 0$ from the theorem 4 we obtain the corresponding results of the works [9], and also [10], if we  take coefficient $\rho \left( t \right)$ at a constant term in the equation as unit.

\newpage

{\large\bf References}

\begin{enumerate}

\item
Kato T. {\it Perturbation theory for linear operators.}
Springer-Verlag, Berlin; Heidelberg, New York, 1966; Mir, Moscow, 1972.

\item
Gasymov M. G., Mirzoev S. S.  {\it On the solvability of the boundary-value problems for the operator-differential equations of elliptic type of the second order.} Differentsial'nye Uravneniya [Differential Equations], {\bf 28} (1992), no. 4, 651-661.

\item
Ilyin V. A., Filippov A. F. {\it About character of a spectrum of self-adjoined extension of Laplace operator in the bounded  area.} Doklady Akad. Nauk SSSR [Soviet Math. Dokl.], {\bf 191} (1970), no. 2, 267-269.

\item
Gorbachuk M. L. {\it Completeness of the system of eigenfunctions and associated functi-ons of a nonself-adjoint boundary value problem for a differential-operator equation of second order.} Funktsional'nyi Analiz I Ego Prilozheniya [Functional Analysis and Its Applications], {\bf 7} (1973), no. 1, 68-69.

\item
Rofe-Beketov F. S. {\it Expansion in eigenfunctions of infinite systems of differential equations in the non-self-adjoint and self-adjoint cases.} Matematicheskii Sbornik [Mathematics of the USSR-Sbornik], {\bf 51(93)} (1960), no. 3, 293-342.

\item
Yakubov S.Y., Aliev B.A. {\it Fredholm property of a boundary value problem with an operator in boundary conditions for an elliptic type operator-differential equation of the second order.} Doklady Akad. Nauk SSSR [Soviet Math. Dokl.], {\bf 257} (1981), no. 5, 1071-1074.

\item
Mirzoyev S.S., Yaqubova  Kh. V. {\it On the solvability of boundary value problems with operators in boundary conditions for one class of operator-differential equations of the second order.} Reports of NAS of Azerbaijan, {\bf 57} (2001), no. 1-3, 12-17.

\item
Aliyev A. R. {\it To the theory of solvability of the second order operator-differential equations with discontinuous coefficients.} Transactions of NAS of Azerb., ser. of phys.-tech. and math. sciences, {\bf 24} (2004), no. 1, 37-44.

\item
Mirzoev S. S. {\it Conditions for the well-defined solvability of boundary-value problems for operator differential equations.} Doklady Akad. Nauk SSSR [Soviet Math. Dokl.], {\bf 273} (1983), no. 2, 292-295.

\item
Aliyev A. R. {\it On the solvability of the boundary-value problem for the operator-differential equations of the third order with discontinuous coefficient.} Proceedings of the Institute Mathematics and Mechanics AS Azerbaijan, {\bf 7(15)} (1997), 18-25.

\item
Aliev A. R. {\it Solubility of boundary-value problems for a class of third-order operator-differential equations in a weighted space.} Uspekhi Matematicheskikh Nauk [Russian Mathematical Surveys], {\bf 60} (2005), no. 4(364), 215-216.

\item
Mirzoyev S. S., Aliyev A. R. {\it Initial boundary value problems for a class of third order operator-differential equations with variable coefficients.} Transactions of NAS of Azerb., ser. of phys.-tech. and math. sciences, {\bf 26} (2006), no. 4, 153-164.

\item
Aliev A. R. {\it On the boundary value problem for a class of operator-differential equations of odd order with variable coefficients.} Doklady Akad. Nauk [Doklady Mathematics], {\bf 421} (2008), no. 2, 151-153.

\item
Lions J. L., Magenes E. {\it Non-homogeneous boundary value problems and applications.} Dunod, Paris, 1968; Mir, Moscow, 1971; Springer-Verlag, Berlin, 1972.

\item
Gorbachuk V. I., Gorbachuk M. L. {\it Boundary value problems for operator differential equations.} Naukova dumka, Kiyev, 1984; Springer, 1990.

\item
Kalyabin G. A. {\it Some problems for Sobolev spaces on the half-line.} Trudy Matematiche-skogo Instituta imeni V. A. Steklova [Proceedings of the Steklov Institute of Mathema-tics], {\bf 255} (2006), 161-169.
\end{enumerate}

\newpage

\begin{titlepage}

{\bf Araz R. Aliev} - the doctor of physical and mathematical sciences

{\bf 1.} Baku State University,

Applied Mathematics and Cybernetics Faculty

23, Z.Khalilov str., AZ1148, Baku, Azerbaijan

{\bf 2.} Institute of Mathematics and Mechanics of NAS of Azerbaijan

9, F.Agayev str., AZ1141, Baku, Azerbaijan

{\bf E-mail:} alievaraz@yahoo.com

{\bf Tel.:} (+994 50) 329 67 38.

{\bf Sevindj F. Babayeva}

Institute of Cybernetics of NAS of Azerbaijan

9, F.Agayev str., AZ1141, Baku, Azerbaijan

{\bf E-mail:} seva\_babaeva@mail.ru

{\bf Tel.:} (+994 50) 348 16 84.

\end{titlepage}

\end{document}